

\input amstex
\expandafter\ifx\csname mathdefs.tex\endcsname\relax
  \expandafter\gdef\csname mathdefs.tex\endcsname{}
\else \message{Hey!  Apparently you were trying to
  \string twice.   This does not make sense.} 
\errmessage{Please edit your file (probably \jobname.tex) and remove
any duplicate ``\string\input'' lines} \fi




\catcode`\X=12\catcode`\@=11

\def\n@wcount{\alloc@0\count\countdef\insc@unt}
\def\n@wwrite{\alloc@7\write\chardef\sixt@@n}
\def\n@wread{\alloc@6\read\chardef\sixt@@n}
\def\r@s@t{\relax}\def\v@idline{\par}\def\@mputate#1/{#1}
\def\l@c@l#1X{\firstpart.#1}\def\gl@b@l#1X{#1}\def\t@d@l#1X{{}}

\def\crossrefs#1{\ifx\all#1\let\tr@ce=\all\else\def\tr@ce{#1,}\fi
   \n@wwrite\cit@tionsout\openout\cit@tionsout=\jobname.cit 
   \write\cit@tionsout{\tr@ce}\expandafter\setfl@gs\tr@ce,}
\def\setfl@gs#1,{\def\@{#1}\ifx\@\empty\let\next=\relax
   \else\let\next=\setfl@gs\expandafter\xdef
   \csname#1tr@cetrue\endcsname{}\fi\next}
\def\m@ketag#1#2{\expandafter\n@wcount\csname#2tagno\endcsname
     \csname#2tagno\endcsname=0\let\tail=\all\xdef\all{\tail#2,}
   \ifx#1\l@c@l\let\tail=\r@s@t\xdef\r@s@t{\csname#2tagno\endcsname=0\tail}\fi
   \expandafter\gdef\csname#2cite\endcsname##1{\expandafter
     \ifx\csname#2tag##1\endcsname\relax?\else\csname#2tag##1\endcsname\fi
     \expandafter\ifx\csname#2tr@cetrue\endcsname\relax\else
     \write\cit@tionsout{#2tag ##1 cited on page \folio.}\fi}
   \expandafter\gdef\csname#2page\endcsname##1{\expandafter
     \ifx\csname#2page##1\endcsname\relax?\else\csname#2page##1\endcsname\fi
     \expandafter\ifx\csname#2tr@cetrue\endcsname\relax\else
     \write\cit@tionsout{#2tag ##1 cited on page \folio.}\fi}
   \expandafter\gdef\csname#2tag\endcsname##1{\expandafter
      \ifx\csname#2check##1\endcsname\relax
      \expandafter\xdef\csname#2check##1\endcsname{}%
      \else\immediate\write16{Warning: #2tag ##1 used more than once.}\fi
      \multit@g{#1}{#2}##1/X%
      \write\t@gsout{#2tag ##1 assigned number \csname#2tag##1\endcsname\space
      on page \number\count0.}%
   \csname#2tag##1\endcsname}}
\def\multit@g#1#2#3/#4X{\def\t@mp{#4}\ifx\t@mp\empty%
      \global\advance\csname#2tagno\endcsname by 1 
      \expandafter\xdef\csname#2tag#3\endcsname
      {#1\number\csname#2tagno\endcsnameX}%
   \else\expandafter\ifx\csname#2last#3\endcsname\relax
      \expandafter\n@wcount\csname#2last#3\endcsname
      \global\advance\csname#2tagno\endcsname by 1 
      \expandafter\xdef\csname#2tag#3\endcsname
      {#1\number\csname#2tagno\endcsnameX}
      \write\t@gsout{#2tag #3 assigned number \csname#2tag#3\endcsname\space
      on page \number\count0.}\fi
   \global\advance\csname#2last#3\endcsname by 1
   \def\t@mp{\expandafter\xdef\csname#2tag#3/}%
   \expandafter\t@mp\@mputate#4\endcsname
   {\csname#2tag#3\endcsname\lastpart{\csname#2last#3\endcsname}}\fi}
\def\t@gs#1{\def\all{}\m@ketag#1e\m@ketag#1s\m@ketag\t@d@l p
   \m@ketag\gl@b@l r \n@wread\t@gsin
   \openin\t@gsin=\jobname.tgs \re@der \closein\t@gsin
   \n@wwrite\t@gsout\openout\t@gsout=\jobname.tgs }
\outer\def\localtags{\t@gs\l@c@l}
\outer\def\globaltags{\t@gs\gl@b@l}
\outer\def\newlocaltag#1{\m@ketag\l@c@l{#1}}
\outer\def\newglobaltag#1{\m@ketag\gl@b@l{#1}}

\newif\ifpr@ 
\def\m@kecs #1tag #2 assigned number #3 on page #4.%
   {\expandafter\gdef\csname#1tag#2\endcsname{#3}
   \expandafter\gdef\csname#1page#2\endcsname{#4}
   \ifpr@\expandafter\xdef\csname#1check#2\endcsname{}\fi}
\def\re@der{\ifeof\t@gsin\let\next=\relax\else
   \read\t@gsin to\t@gline\ifx\t@gline\v@idline\else
   \expandafter\m@kecs \t@gline\fi\let \next=\re@der\fi\next}
\def\pretags#1{\pr@true\pret@gs#1,,}
\def\pret@gs#1,{\def\@{#1}\ifx\@\empty\let\n@xtfile=\relax
   \else\let\n@xtfile=\pret@gs \openin\t@gsin=#1.tgs \message{#1} \re@der 
   \closein\t@gsin\fi \n@xtfile}

\newcount\sectno\sectno=0\newcount\subsectno\subsectno=0
\newif\ifultr@local \def\ultralocal{\ultr@localtrue}
\def\firstpart{\number\sectno}
\def\lastpart#1{\ifcase#1 \or a\or b\or c\or d\or e\or f\or g\or h\or 
   i\or k\or l\or m\or n\or o\or p\or q\or r\or s\or t\or u\or v\or w\or 
   x\or y\or z \fi}

\def\resetall{\global\advance\sectno by 1\subsectno=0
   \gdef\firstpart{\number\sectno}\r@s@t}
\def\resetsub{\global\advance\subsectno by 1
   \gdef\firstpart{\number\sectno.\number\subsectno}\r@s@t}
\def\newsection#1\par{\resetall\vskip0pt plus.3\vsize\penalty-250
   \vskip0pt plus-.3\vsize\bigskip\bigskip
   \message{#1}\leftline{\bf#1}\nobreak\bigskip}
\def\subsection#1\par{\ifultr@local\resetsub\fi
   \vskip0pt plus.2\vsize\penalty-250\vskip0pt plus-.2\vsize
   \bigskip\smallskip\message{#1}\leftline{\bf#1}\nobreak\medskip}

\def\t@gsoff#1,{\def\@{#1}\ifx\@\empty\let\next=\relax\else\let\next=\t@gsoff
   \def\@@{p}\ifx\@\@@\else
   \expandafter\gdef\csname#1cite\endcsname##1{\zeigen{##1}}
   \expandafter\gdef\csname#1page\endcsname##1{?}
   \expandafter\gdef\csname#1tag\endcsname##1{\zeigen{##1}}\fi\fi\next}
\def\verbatimtags{\ifx\all\relax\else\expandafter\t@gsoff\all,\fi}
\def\zeigen#1{\hbox{$\langle$}#1\hbox{$\rangle$}}

\def\(#1){\edef\dot@g{\ifmmode\ifinner(\hbox{\noexpand\etag{#1}})
   \else\noexpand\eqno(\hbox{\noexpand\etag{#1}})\fi
   \else(\noexpand\ecite{#1})\fi}\dot@g}

\newif\ifbr@ck
\def\eat#1{}
\def\[#1]{\br@cktrue[\br@cket#1'X]}
\def\br@cket#1'#2X{\def\temp{#2}\ifx\temp\empty\let\next\eat
   \else\let\next\br@cket\fi
   \ifbr@ck\br@ckfalse\br@ck@t#1,X\else\br@cktrue#1\fi\next#2X}
\def\br@ck@t#1,#2X{\def\temp{#2}\ifx\temp\empty\let\neext\eat
   \else\let\neext\br@ck@t\def\temp{,}\fi
   \def\teemp{#1}\ifx\teemp\empty\else\rcite{#1}\fi\temp\neext#2X}
\def\resetbr@cket{\gdef\[##1]{[\rtag{##1}]}}
\def\references{\resetbr@cket\newsection References\par}

\newtoks\symb@ls\newtoks\s@mb@ls\newtoks\p@gelist\n@wcount\ftn@mber
    \ftn@mber=1\newif\ifftn@mbers\ftn@mbersfalse\newif\ifbyp@ge\byp@gefalse
\def\defm@rk{\ifftn@mbers\n@mberm@rk\else\symb@lm@rk\fi}
\def\n@mberm@rk{\xdef\m@rk{{\the\ftn@mber}}%
    \global\advance\ftn@mber by 1 }
\def\rot@te#1{\let\temp=#1\global#1=\expandafter\r@t@te\the\temp,X}
\def\r@t@te#1,#2X{{#2#1}\xdef\m@rk{{#1}}}
\def\b@@st#1{{$^{#1}$}}\def\str@p#1{#1}
\def\symb@lm@rk{\ifbyp@ge\rot@te\p@gelist\ifnum\expandafter\str@p\m@rk=1 
    \s@mb@ls=\symb@ls\fi\write\f@nsout{\number\count0}\fi \rot@te\s@mb@ls}
\def\byp@ge{\byp@getrue\n@wwrite\f@nsin\openin\f@nsin=\jobname.fns 
    \n@wcount\currentp@ge\currentp@ge=0\p@gelist={0}
    \re@dfns\closein\f@nsin\rot@te\p@gelist
    \n@wread\f@nsout\openout\f@nsout=\jobname.fns }
\def\m@kelist#1X#2{{#1,#2}}
\def\re@dfns{\ifeof\f@nsin\let\next=\relax\else\read\f@nsin to \f@nline
    \ifx\f@nline\v@idline\else\let\t@mplist=\p@gelist
    \ifnum\currentp@ge=\f@nline
    \global\p@gelist=\expandafter\m@kelist\the\t@mplistX0
    \else\currentp@ge=\f@nline
    \global\p@gelist=\expandafter\m@kelist\the\t@mplistX1\fi\fi
    \let\next=\re@dfns\fi\next}
\def\symbols#1{\symb@ls={#1}\s@mb@ls=\symb@ls} 
\def\bigsymbol{\textstyle}
\symbols{\bigsymbol\ast,\dagger,\ddagger,\sharp,\flat,\natural,\star}
\def\ftnumbers{\ftn@mberstrue} \def\ftsymbols{\ftn@mbersfalse}
\def\paginal{\byp@ge} \def\resetftnumbers{\ftn@mber=1}
\def\ftnote#1{\defm@rk\expandafter\expandafter\expandafter\footnote
    \expandafter\b@@st\m@rk{#1}}

\long\def\jump#1\endjump{}
\def\ssum{\mathop{\lower .1em\hbox{$\textstyle\Sigma$}}\nolimits}

\def\qed{\nobreak\kern 1em \vrule height .5em width .5em depth 0em}
\def\newneq{\hbox{\rlap{\hbox to 1\wd9{\hss$=$\hss}}\raise .1em 
   \hbox to 1\wd9{\hss$\scriptscriptstyle/$\hss}}}
\def\subsetne{\setbox9 = \hbox{$\subset$}\mathrel{\hbox{\rlap
   {\lower .4em \newneq}\raise .13em \hbox{$\subset$}}}}
\def\supsetne{\setbox9 = \hbox{$\subset$}\mathrel{\hbox{\rlap
   {\lower .4em \newneq}\raise .13em \hbox{$\supset$}}}}

\def\vbar{\mathchoice{\vrule height6.3ptdepth-.5ptwidth.8pt\kern-.8pt}
   {\vrule height6.3ptdepth-.5ptwidth.8pt\kern-.8pt}
   {\vrule height4.1ptdepth-.35ptwidth.6pt\kern-.6pt}
   {\vrule height3.1ptdepth-.25ptwidth.5pt\kern-.5pt}}
\def\f@dge{\mathchoice{}{}{\mkern.5mu}{\mkern.8mu}}
\def\b@c#1#2{{\rm \mkern#2mu\vbar\mkern-#2mu#1}}
\def\b@b#1{{\rm I\mkern-3.5mu #1}}
\def\b@a#1#2{{\rm #1\mkern-#2mu\f@dge #1}}
\def\bb#1{{\count4=`#1 \advance\count4by-64 \ifcase\count4\or\b@a A{11.5}\or
   \b@b B\or\b@c C{5}\or\b@b D\or\b@b E\or\b@b F \or\b@c G{5}\or\b@b H\or
   \b@b I\or\b@c J{3}\or\b@b K\or\b@b L \or\b@b M\or\b@b N\or\b@c O{5} \or
   \b@b P\or\b@c Q{5}\or\b@b R\or\b@a S{8}\or\b@a T{10.5}\or\b@c U{5}\or
   \b@a V{12}\or\b@a W{16.5}\or\b@a X{11}\or\b@a Y{11.7}\or\b@a Z{7.5}\fi}}

\catcode`\X=11 \catcode`\@=12

\expandafter\ifx\csname citeadd.tex\endcsname\relax
\expandafter\gdef\csname citeadd.tex\endcsname{}
\else \message{Hey!  Apparently you were trying to
\string twice.   This does not make sense.} 
\errmessage{Please edit your file (probably \jobname.tex) and remove
any duplicate ``\string\input'' lines} \fi

\sectno=-1   
\localtags
\def\cite #1{\rm[#1]}
\NoBlackBoxes
\documentstyle {amsppt}
\topmatter
\title {Erd\"os and R\'enyi Conjecture} \endtitle
\author {Saharon Shelah \thanks{\null\newline
Latest Revision 97/Aug/14 \newline
I thank Alice Leonhardt for the beautiful typing} \endthanks} \endauthor
\affil{Institute of Mathematics \\
The Hebrew University \\
Jerusalem, Israel
\medskip
Rutgers University \\
Department of Mathematics \\
New Brunswick, NJ  USA} \endaffil
\abstract{Affirming a conjecture of Erd\"os and R\'enyi we prove that for
any (real number) $c_1 > 0$ for some $c_2 > 0$, if a graph $G$ has no 
$c_1(\text{log } n)$ nodes on which the graph is complete or edgeless (i.e.
$G$ exemplifies $|G| \nrightarrow (c_1 \text{ log } n)^2_2$) \underbar{then}
$G$ has at least $2^{c_2n}$ non-isomorphic (induced) subgraphs.} \endabstract
\endtopmatter
\document  

\expandafter\ifx\csname alice2jlem.tex\endcsname\relax
  \expandafter\gdef\csname alice2jlem.tex\endcsname{}
\else \message{Hey!  Apparently you were trying to
\string  twice.   This does not make sense.}
\errmessage{Please edit your file (probably \jobname.tex) and remove
any duplicate ``\string\input'' lines} \fi

\expandafter\ifx\csname bib4plain.tex\endcsname\relax
  \expandafter\gdef\csname bib4plain.tex\endcsname{}
\else \message{Hey!  Apparently you were trying to \string twice.   This does not make sense.}
\errmessage{Please edit your file (probably \jobname.tex) and remove
any duplicate ``\string\input'' lines} \fi

\def\renewcommand{\newcommand}	       
\edef\cite{\the\catcode`@}%
\catcode`@ = 11
\let\@oldatcatcode = \cite
\chardef\@letter = 11
\chardef\@other = 12
%
%
%
%
\def\@innerdef#1#2{\edef#1{\expandafter\noexpand\csname #2\endcsname}}%
%
%
\@innerdef\@innernewcount{newcount}%
\@innerdef\@innernewdimen{newdimen}%
\@innerdef\@innernewif{newif}%
\@innerdef\@innernewwrite{newwrite}%
%
%
%
\def\@gobble#1{}%
%
%
%
\ifx\inputlineno\@undefined
   \let\@linenumber = \empty 
\else
   \def\@linenumber{\the\inputlineno:\space}%
\fi
%
%
%
\def\@futurenonspacelet#1{\def\cs{#1}%
   \afterassignment\@stepone\let\@nexttoken=
}%
\begingroup 
\def\\{\global\let\@stoken= }%
\\ 
\endgroup
\def\@stepone{\expandafter\futurelet\cs\@steptwo}%
\def\@steptwo{\expandafter\ifx\cs\@stoken\let\@@next=\@stepthree
   \else\let\@@next=\@nexttoken\fi \@@next}%
\def\@stepthree{\afterassignment\@stepone\let\@@next= }%
%
%
%
\def\@getoptionalarg#1{%
   \let\@optionaltemp = #1%
   \let\@optionalnext = \relax
   \@futurenonspacelet\@optionalnext\@bracketcheck
}%
%
%
\def\@bracketcheck{%
   \ifx [\@optionalnext
      \expandafter\@@getoptionalarg
   \else
      \let\@optionalarg = \empty
      \expandafter\@optionaltemp
   \fi
}%
\def\@@getoptionalarg[#1]{%
   \def\@optionalarg{#1}%
   \@optionaltemp
}%
%
%
%
\def\@nnil{\@nil}%
\def\@fornoop#1\@@#2#3{}%
\def\@for#1:=#2\do#3{%
   \edef\@fortmp{#2}%
   \ifx\@fortmp\empty \else
      \expandafter\@forloop#2,\@nil,\@nil\@@#1{#3}%
   \fi
}%
\def\@forloop#1,#2,#3\@@#4#5{\def#4{#1}\ifx #4\@nnil \else
       #5\def#4{#2}\ifx #4\@nnil \else#5\@iforloop #3\@@#4{#5}\fi\fi
}%
\def\@iforloop#1,#2\@@#3#4{\def#3{#1}\ifx #3\@nnil
       \let\@nextwhile=\@fornoop \else
      #4\relax\let\@nextwhile=\@iforloop\fi\@nextwhile#2\@@#3{#4}%
}%
%
%
%
\@innernewif\if@fileexists
\def\@testfileexistence{\@getoptionalarg\@finishtestfileexistence}%
\def\@finishtestfileexistence#1{%
   \begingroup
      \def\extension{#1}%
      \immediate\openin0 =
         \ifx\@optionalarg\empty\jobname\else\@optionalarg\fi
         \ifx\extension\empty \else .#1\fi
         \space
      \ifeof 0
         \global\@fileexistsfalse
      \else
         \global\@fileexiststrue
      \fi
      \immediate\closein0
   \endgroup
}%
%
%
%
%
\def\bibliographystyle#1{%
   \@readauxfile
   \@writeaux{\string\bibstyle{#1}}%
}%
\let\bibstyle = \@gobble
%
%
\let\bblfilebasename = \jobname
\def\bibliography#1{%
   \@readauxfile
   \@writeaux{\string\bibdata{#1}}%
   \@testfileexistence[\bblfilebasename]{bbl}%
   \if@fileexists
      \nobreak
      \@readbblfile
   \fi
}%
\let\bibdata = \@gobble
%
%
\def\nocite#1{%
   \@readauxfile
   \@writeaux{\string\citation{#1}}%
}%
\@innernewif\if@notfirstcitation
%
%
\def\cite{\@getoptionalarg\@cite}%
%
%
\def\@cite#1{%
   \let\@citenotetext = \@optionalarg
   \printcitestart
   \nocite{#1}%
   \@notfirstcitationfalse
   \@for \@citation :=#1\do
   {%
      \expandafter\@onecitation\@citation\@@
   }%
   \ifx\empty\@citenotetext\else
      \printcitenote{\@citenotetext}%
   \fi
   \printcitefinish
}%
\def\@onecitation#1\@@{%
   \if@notfirstcitation
      \printbetweencitations
   \fi
   \expandafter \ifx \csname\@citelabel{#1}\endcsname \relax
      \if@citewarning
         \message{\@linenumber Undefined citation `#1'.}%
      \fi
      \expandafter\gdef\csname\@citelabel{#1}\endcsname{%
\strut
\vadjust{\vskip-\dp\strutbox
\vbox to 0pt{\vss\parindent0cm \leftskip=\hsize 
\advance\leftskip3mm
\advance\hsize 4cm\strut\openup-4pt 
\rightskip 0cm plus 1cm minus 0.5cm ?  #1 ?\strut}}
         {\tt
            \escapechar = -1
            \nobreak\hskip0pt
            \expandafter\string\csname#1\endcsname
            \nobreak\hskip0pt
         }%
      }%
   \fi
   \csname\@citelabel{#1}\endcsname
   \@notfirstcitationtrue
}%
%
%
\def\@citelabel#1{b@#1}%
%
%
\def\@citedef#1#2{\expandafter\gdef\csname\@citelabel{#1}\endcsname{#2}}%
%
%
%
\def\@readbblfile{%
   \ifx\@itemnum\@undefined
      \@innernewcount\@itemnum
   \fi
   \begingroup
      \def\begin##1##2{%
         \setbox0 = \hbox{\biblabelcontents{##2}}%
         \biblabelwidth = \wd0
      }%
      \def\end##1{}
      %
      %
      \@itemnum = 0
      \def\bibitem{\@getoptionalarg\@bibitem}%
      \def\@bibitem{%
         \ifx\@optionalarg\empty
            \expandafter\@numberedbibitem
         \else
            \expandafter\@alphabibitem
         \fi
      }%
      \def\@alphabibitem##1{%
         \expandafter \xdef\csname\@citelabel{##1}\endcsname {\@optionalarg}%
         \ifx\biblabelprecontents\@undefined
            \let\biblabelprecontents = \relax
         \fi
         \ifx\biblabelpostcontents\@undefined
            \let\biblabelpostcontents = \hss
         \fi
         \@finishbibitem{##1}%
      }%
      \def\@numberedbibitem##1{%
         \advance\@itemnum by 1
         \expandafter \xdef\csname\@citelabel{##1}\endcsname{\number\@itemnum}%
         \ifx\biblabelprecontents\@undefined
            \let\biblabelprecontents = \hss
         \fi
         \ifx\biblabelpostcontents\@undefined
            \let\biblabelpostcontents = \relax
         \fi
         \@finishbibitem{##1}%
      }%
      \def\@finishbibitem##1{%
         \biblabelprint{\csname\@citelabel{##1}\endcsname}%
         \@writeaux{\string\@citedef{##1}{\csname\@citelabel{##1}\endcsname}}%
         \ignorespaces
      }%
      %
      %
      \let\em = \bblem
      \let\newblock = \bblnewblock
      \let\sc = \bblsc
      \frenchspacing
      \clubpenalty = 4000 \widowpenalty = 4000
      \tolerance = 10000 \hfuzz = .5pt
      \everypar = {\hangindent = \biblabelwidth
                      \advance\hangindent by \biblabelextraspace}%
      \bblrm
      \parskip = 1.5ex plus .5ex minus .5ex
      \biblabelextraspace = .5em
      \bblhook
      \input \bblfilebasename.bbl
   \endgroup
}%
%
%
\@innernewdimen\biblabelwidth
\@innernewdimen\biblabelextraspace
%
%
%
\def\biblabelprint#1{%
   \noindent
   \hbox to \biblabelwidth{%
      \biblabelprecontents
      \biblabelcontents{#1}%
      \biblabelpostcontents
   }%
   \kern\biblabelextraspace
}%
%
%
%
\def\biblabelcontents#1{{\bblrm [#1]}}%
%
%
\def\bblrm{\rm}%
%
%
\def\bblem{\it}%
%
%
\def\bblsc{\ifx\@scfont\@undefined
              \font\@scfont = cmcsc10
           \fi
           \@scfont
}%
%
%
\def\bblnewblock{\hskip .11em plus .33em minus .07em }%
%
%
\let\bblhook = \empty
%
%
%
\def\printcitestart{[}
\def\printcitefinish{]}
\def\printbetweencitations{, }
\def\printcitenote#1{, #1}
%
%
%
\let\citation = \@gobble
%
%
%
\@innernewcount\@numparams
%
%
\def\newcommand#1{%
   \def\@commandname{#1}%
   \@getoptionalarg\@continuenewcommand
}%
%
%
\def\@continuenewcommand{%
   \@numparams = \ifx\@optionalarg\empty 0\else\@optionalarg \fi \relax
   \@newcommand
}%
%
%
\def\@newcommand#1{%
   \def\@startdef{\expandafter\edef\@commandname}%
   \ifnum\@numparams=0
      \let\@paramdef = \empty
   \else
      \ifnum\@numparams>9
         \errmessage{\the\@numparams\space is too many parameters}%
      \else
         \ifnum\@numparams<0
            \errmessage{\the\@numparams\space is too few parameters}%
         \else
            \edef\@paramdef{%
               \ifcase\@numparams
                  \empty  No arguments.
               \or ####1%
               \or ####1####2%
               \or ####1####2####3%
               \or ####1####2####3####4%
               \or ####1####2####3####4####5%
               \or ####1####2####3####4####5####6%
               \or ####1####2####3####4####5####6####7%
               \or ####1####2####3####4####5####6####7####8%
               \or ####1####2####3####4####5####6####7####8####9%
               \fi
            }%
         \fi
      \fi
   \fi
   \expandafter\@startdef\@paramdef{#1}%
}%
%
%
%
%
\def\@readauxfile{%
   \if@auxfiledone \else 
      \global\@auxfiledonetrue
      \@testfileexistence{aux}%
      \if@fileexists
         \begingroup
            \endlinechar = -1
            \catcode`@ = 11
            \input \jobname.aux
         \endgroup
      \else
         \message{\@undefinedmessage}%
         \global\@citewarningfalse
      \fi
      \immediate\openout\@auxfile = \jobname.aux
   \fi
}%
%
%
\newif\if@auxfiledone
\ifx\noauxfile\@undefined \else \@auxfiledonetrue\fi
%
%
%
%
\@innernewwrite\@auxfile
\def\@writeaux#1{\ifx\noauxfile\@undefined \write\@auxfile{#1}\fi}%
%
%
%
\ifx\@undefinedmessage\@undefined
   \def\@undefinedmessage{No .aux file; I won't give you warnings about
                          undefined citations.}%
\fi
%
%
\@innernewif\if@citewarning
\ifx\noauxfile\@undefined \@citewarningtrue\fi
%
%
%
\catcode`@ = \@oldatcatcode


\def\widestnumber#1#2{}

\def\rm{\fam0 \tenrm}

\def\fakesubhead#1\endsubhead{\bigskip\noindent{\bf#1}\par}


%
%
%

%

\font\textrsfs=rsfs10
\font\scriptrsfs=rsfs7
\font\scriptscriptrsfs=rsfs5

\newfam\rsfsfam
\textfont\rsfsfam=\textrsfs
\scriptfont\rsfsfam=\scriptrsfs
\scriptscriptfont\rsfsfam=\scriptscriptrsfs

\edef\oldcatcodeofat{\the\catcode`\@}
\catcode`\@11

\def\Cal@@#1{\noaccents@ \fam \rsfsfam #1}

\catcode`\@\oldcatcodeofat

\newpage

\head {\S0 Introduction} \endhead  \resetall 
\bigskip

Erd\"os and R\'enyi conjectured (letting $I(G)$ denote the number 
of (induced) subgraphs of $G$ up to isomorphism and $Rm(G)$ be the 
maximal number of nodes on which $G$ is complete or edgeless):
\medskip
\roster
\item "{$(*)$}" for every $c_1 > 0$ for some $c_2 > 0$ for $n$ large enough
for every graph $G_n$ with $n$ points
{\roster  
\itemitem{ $\bigotimes$ }  $Rm(G_n) < c_1(\text{log } n) \Rightarrow I(G_n)
\ge 2^{c_2n}$.
\endroster}
\endroster
\medskip

\noindent
They succeeded to prove a parallel theorem replacing $Rm(G)$ by the bipartite 
version:

$$
\align
\text{Bipartite}(G) =: \text{ Max} \biggl\{ k:&\text{ there are disjoint 
sets } A_1,A_2 \text{ of } k \text{ nodes of } G, \\
  &\text{ such that } (\forall x_1 \in A_1)(\forall x_2 \in A_2)
(\{x_1,x_2\} \text{ an edge) or} \\
  &\,(\forall x_1 \in A_1)(\forall x_2 \in A_2)(\{x_1,x_2\} \text{ is not an
edge)} \biggr\}.
\endalign
$$ 
\medskip

\noindent
It is well known that $Rm(G_n) \ge \frac12 \text{ log } n$.  On the other
hand, Erd\"os \cite{Er7} proved that for every $n$ for some graph 
$G_n,Rm(G_n) \le 2 \text{ log } n$.  In his construction $G_n$ is quite a
random graph; it seems reasonable that any graph $G_n$ with 
small Rm$(G_n)$ is of similar character and this is the rationale of the
conjecture.

Alon and Bollobas \cite{AlBl} and Erd\"os and Hajnal \cite{EH9} affirm a
conjecture of Hajnal:
\medskip
\roster
\item "{$(*)$}"  if $Rm(G_n) < (1-\varepsilon)n$ then $I(G_n) > \Omega
(\varepsilon n^2)$ \newline
and Erd\"os and Hajnal \cite{EH9} also prove
\item "{$(*)$}"  for any fixed $k$, if $Rm(G_n) < \frac nk$ then
$I(G_n) > n^{\Omega(\sqrt k)}$.
\endroster
\medskip

\noindent
Alon and Hajnal \cite{AH} noted that those results give poor bounds for
$I(G_n)$ in the case Rm$(G_n)$ is much smaller than a multiple of log $n$, and
prove an inequality weaker than the conjecture:

$$
I(G_n) \ge 2^{n/2t^{20\text{ log}(2t)}} 
\text{ when } t = Rm(G_m) \tag"{$(*)$}"
$$
\medskip

\noindent
so in particular if $t \ge c \text{ log } n$ they got $I(G_n) \ge 
2^{n/(\text{log }n)^{c \text{ log log }n}}$, that is the constant $c_2$ 
in the conjecture is replaced by (log $n)^{c \text{ log log } n}$ for 
some $c$.

I thank Andras Hajnal for telling me about the problem and 
Mariusz Rabus and Andres Villaveces for some corrections.
\newpage

\head {\S1} \endhead  \resetall
\bigskip

\demo{\stag{0.1}  Notation}  log $n = \text{ log}_2n$. \newline
Let $c$ denote a positive real. \newline
$G,H$ denote graphs, which are here finite, simple and undirected. \newline
$V^G$ is the set of nodes of the graph $G$. \newline
$E^G$ is the set of edges of the graph $G$ so $G = (V^G,E^G),E^G$ is a
symmetric, irreflexive relation on $V^G$ i.e. a set of unordered pairs.  So
$\{x,y\} \in E^G,xEy,\{x,y\}$ an edge of $G$, all have the same meaning.
\newline
$H \subseteq G$ means
that $H$ is an induced subgraph of $G$; i.e. 
$H = G \restriction V^H$. \newline  
Let $|X|$ be the number of elements of the set $X$.
\enddemo
\bigskip

\definition{\stag{1} Definition}  $I(G)$ is the number of (induced)
subgraphs of $G$ up to isomorphisms.
\enddefinition
\bigskip

\proclaim{\stag{2} Theorem}  For any $c_1 > 0$ for some $c_2 > 0$ we have
(for $n$ large enough): if $G$ is a graph with $n$ edges and $G$ has neither
a complete subgraph with $\ge c_1 \text{ log }n$ nodes nor a subgraph with
no edges with $\ge c_1 \text{ log } n$ nodes \underbar{then} $I(G) \ge
2^{c_2n}$.
\endproclaim
\bigskip

\remark{\stag{3} Remark}  1) Suppose $n \nrightarrow (r_1,r_2)$ and 
$m$ are given.  Choose a graph $H$ on \newline
$\{0,\dotsc,n-1\}$ exemplifying 
$n \nrightarrow (r_1,r_2)^2$ (i.e. with no complete subgraphs with $r_1$ nodes
and no independent set with $r_2$ nodes).  Define the graph $G$ with set of
nodes $V^G = \{0,\dotsc,mn-1\}$ and set of edges $E^G = \{\{mi_1 + \ell_1,
mi_2 + \ell_2\}:\{i_1,i_2\} \in E^H$ and $\ell_1,\ell_2 < m\}$.  Clearly
$G$ has $nm$ nodes and it exemplifies $mn \nrightarrow (r_1,mr_2)$.  So
$I(G) \le (m+1)^n \le 2^{n\text{ log}_2(m+1)}$ (as the isomorphism type of
$G' \subseteq G$ is determined by $\langle |G' \cap [mi,mi+m)|:i < n
\rangle$).  We conjecture that this is the worst case. \newline
2) Similarly if $n \nrightarrow \left( \bmatrix r_1 \\ r_2 
\endbmatrix \right)^2_2$; i.e.
there is a graph with $n$ nodes and no disjoint $A_1,A_2 \subseteq V^G,
|A_1| = r_1,|A_2| = r_2$ such that $A_1 \times A_2 \subseteq E^G$ or 
$(A_1 \times A_2) \cap E^G = \emptyset$, \underbar{then} there is $G$ 
exemplifying $mn \rightarrow \left( \bmatrix n_1 m \\ r_2 m \endbmatrix 
\right)^2_2$ such that $I(G) \le 2^{n\text{ log}(m+1)}$.
\endremark
\bigskip

\demo{Proof}  Let $c_1$, a real $> 0$, be given. \newline
Let $m^*_1$ be \footnote{the log log $n$ can be replaced by a constant 
computed from $m^*_1,m^*_2,c_\ell$ later} such that for every $n$ (large
enough)
${\frac n{(\text{log } n)^2\text{log log}n}} \rightarrow (c_1 \text{ log }
n,\frac {c_1}{m^*_1} \text{ log } n)$. \newline
[Why does it exist?  
By Erd\"os and Szekeres \cite{ErSz} $\binom {n_1+n_2-2}{n-1} 
\rightarrow (n_1,n_2)^2$ and hence for any $k$ letting $n_1=km,n_2=m$ we have
$\binom{km+m-2}{m-1} \rightarrow (km,m)^2$, now
$\binom{m+m-2}{m-1} \le 2^{2(m-1)}$ and

$$
\align
\binom{(k+1)m+m-2}{m-1} \bigl/ \binom{km+m-2}{m-1} &= 
\dsize \prod^{m-2}_{i=0} (1 + \frac{m}{km+i}) \\
  &\le \dsize \prod^{m-2}_{i=0} (1+ \frac{m}{km}) =
(1 + \frac 1k)^{m-1}
\endalign
$$
\medskip

\noindent
hence $\binom{km+m-2}{m-1} \le \left( 4 \cdot \dsize \prod^{k-2}_{\ell = 0}
(1 + \frac{1}{\ell + 1}) \right)^{m-1}$, and choose $k$ large enough 
(see below).  For (large enough) $n$ we let 
$m = (c_1 \log n)/k$, more exactly the first integer is not below this number
so

$$
\align
\text{log}\binom{km+m-2}{m-1} &\le \text{ log}
\left( 4 \cdot \dsize \prod^{k-2}_{\ell = 0} (1 + \frac{1}{\ell +1}) \right)
^{m-1} \\
  &\le (\text{log } n) \cdot \frac{c_1}{k} \cdot \text{ log}
\left( 4 \cdot \dsize \prod^{k-2}_{\ell =0} (1 + \frac{1}{\ell + 1}) \right)
\le \frac 12 (\text{log } n)
\endalign
$$
\medskip

\noindent
(the last inequality holds as $k$ is large enough); lastly let $m^*_1$ be 
such a $k$. Alternatively, just repeat the proof of Ramsey's theorem.]
\medskip

\noindent
Let $m^*_2$ be minimal such that $m^*_2 \rightarrow (m^*_1)^2_2$.
\medskip

\noindent
Let $c_2 < \frac 1{m^*_2}$ (be a positive real).
\medskip

\noindent
Let $c_3 \in (0,1)_{\Bbb R}$ be such that $0 < c_3 < \frac 1{m^*_2} - c_2$.
\medskip

\noindent
Let $c_4 \in \Bbb R^+$ be $4/c_3$ (even $(2 + \varepsilon)/c_3$ suffices).
\medskip

\noindent
Let $c_5 = \frac {1-c_2-c_3}{m^*_2}$ (it is $> 0$).
\medskip

\noindent
Let $\varepsilon \in (0,1)_{\Bbb R}$ be small enough. 
\medskip
\noindent
Now suppose
\medskip
\roster
\item "{$(*)_0$}"  $n$ is large enough, $G$ a graph with $n$ nodes and
$I(G) < 2^{c_2n}$.
\endroster
\medskip

We choose $A \subseteq V^G$ in the following random way: for each $x \in
V^G$ we flip a coin with probability $c_3/\text{log }n$, and let $A$ be the
set of $x \in V^G$ for which we succeed.  For any $A \subseteq V^G$ let
$\approx_A$ be the following relation on $V^G,x \approx_A y$ \underbar{iff}
$x,y \in V^G$ and $(\forall z \in A)[z E^G x \leftrightarrow z E^Gy]$.  
Clearly $\approx_A$ is an equivalence relation; and let $\approx'_A = 
\approx_A \restriction (V^G \backslash A)$.

For distinct $x,y \in V^G$ what is the probability that $x \approx_A y$?  Let

$$
\text{Dif}(x,y) =: \{z:z \in V^G \text{ and } z E^G x \leftrightarrow \neg
z E^G y\},
$$
\medskip

\noindent
and dif$(x,y) = |\text{Dif}(x,y)|$, so the probability of
$x \approx_A y$ is \newline
$\left( 1 - \frac {c_3}{\text{log }n} \right)^{\text{dif}(x,y)} 
\sim e^{-c_3 \text{ dif}(x,y)/\text{log }n}$. \newline
Hence the probability that for some $x \ne y$ in $V^G$ satisfying 
dif$(x,y) \ge c_4 (\text{log }n)^2$ we have $x \approx_A y$ is at most

$$
\binom n2 e^{-c_3(c_4(\text{log }n)^2)/\text{log }n} \le \binom n2
e^{-4 \text{ log } n} \le 1/n^2
$$
\medskip

\noindent
(remember $c_3c_4 = 4$ and $(4/ \text{log } e) \ge 2$).  Hence for some set
$A$ of nodes of $G$ we have
\medskip
\roster
\item "{$(*)_1$}"  $A \subseteq V^G$ and $A$ has 
$\le \frac {c_3}{\text{log } n} \cdot n$ elements and $A$ is non-empty and
\smallskip
\noindent
\item "{$(*)_2$}"  if $x \approx_A y$ \underbar{then} dif$(x,y) \le c_4
(\text{log }n)^2$.
\endroster
\medskip

\noindent
Next
\medskip
\roster
\item "{$(*)_3$}"  $\ell =: |(V^G \backslash A)/\approx_A|$ (i.e. the number
of equivalence classes of \newline
$\approx'_A = \approx_A \restriction (V^G \backslash
A)$) is $< (c_2 + c_3) \cdot n$ \newline
\smallskip
\noindent
[why?  let $C_1,\dotsc,C_\ell$ be the $\approx'_A$-equivalence classes.  
For each $u \subseteq \{1,\dotsc,\ell\}$ let $G_u = G \restriction (A \cup
\dsize \bigcup_{i \in u}C_i)$.  So $G_u$ is an induced subgraph of $G$ and
$(G_u,c)_{c \in A}$ for $u \subseteq \{1,\dotsc,\ell\}$ are pairwise non-
isomorphic structures, so
\endroster

$$
\align
2^\ell = |\{u:u \subseteq \{1,\dotsc,\ell\}\}| &\le |\{f:f \text{ a function
from } A \text{ into } V^G\}| \times I(G) \\
  &\le n^{|A|} \times I(G),
\endalign
$$
\medskip

\noindent
hence (first inequality by the hypothesis toward contradiction)

$$
\align
2^{c_2n} > I(G) \ge 2^\ell \times n^{-|A|} &\ge 2^\ell \cdot
n^{-c_3n/\text{log }n} \\
  &= 2^\ell \times 2^{-c_3n}
\endalign
$$
\medskip

\noindent
hence 

$$
c_2n > \ell - c_3n \text{ so } \ell < (c_2 + c_3)n \text{ and we have
gotten } (*)_3].
$$
\medskip

\noindent
Let $\{B_i:i < i^*\}$ be a maximal family such that:
\medskip
\roster
\item "{$(a)$}"  each $B_i$ is a subset of some $\approx'_A$-equivalence
class
\smallskip
\item "{$(b)$}"  the $B_i$'s are pairwise disjoint
\smallskip
\item "{$(c)$}"  $|B_i| = m^*_1$
\smallskip
\item "{$(d)$}"  $G \restriction B_i$ is a complete graph or a graph with
no edges.
\endroster
\medskip

\noindent
Now if $x \in V^G \backslash A$ then $(x/\approx'_A) \backslash \dsize
\bigcup_{i < i^*} B_i$ has $< m^*_2$ elements (as $m^*_2 \rightarrow
(m^*_1)^2_2$ by the choice of $m^*_2$ and ``$\langle B_i:i < i^* \rangle$ is
maximal").  Hence

$$
\align
n = |V^G| &= |A| + |\dsize \bigcup_{i < i^*} B_i| + |V^G \backslash A
\backslash \dsize \bigcup_{i < i^*} B_i| \\
  &\le c_3 \frac n{\text{log } n} + m^*_1 \times i^* +
|(V^G \backslash A)/\approx'_A| \times m^*_2 \\
  &\le c_3 \frac n{\text{log } n} + m^*_1 \times i^* + m^*_2(c_2 + c_3)n \\
  &= c_3 \frac n{\text{log } n} + m^*_1 \times i^* + (1 - m^*_2c_5) \cdot n
\endalign
$$
\medskip

\noindent
hence
\medskip
\roster
\item "{$(*)_4$}"  $i^* \ge \frac{n}{m^*_1}(m^*_2c_5 - \frac{c_3}
{\text{log }n})$.
\endroster
\medskip

\noindent
For $i < i^*$ let

$$
B_i = \{x_{i,0},x_{i,2},\dotsc,x_{m^*_1-1}\},
$$
\medskip

\noindent
and let  

$$
\align
u_i =: \biggl\{ j < i^*:&j \ne i \text{ and for some } \ell_1 \in
\{1,\dotsc,m^*_1-1\} \text{ and} \\
  &\ell_2 \in \{0,\dotsc,m^*_1-1\} \text{ we have} \\
  &x_{j,\ell_2} \in \text{ Dif}(x_{i,0},x_{i,\ell_1}) \biggr\}.
\endalign
$$
\medskip

\noindent
Clearly
\medskip
\roster
\item "{$(*)_5$}"  $|u_i| \le m^*_1(m^*_1 - 1)c_4(\text{log } n)^2$.
\endroster
\medskip

\noindent
Next we can find $W$ such that
\medskip

$(*)_6 \quad (i) \,\,W \subseteq \{0,\dotsc,i^*-1\}$

$\quad \quad \,\,(ii) \,\,\,|W| \ge i^*/(m^*_1(m^*_1 - 1)
c_4(\text{log } n)^2)$

$\quad \quad (iii) \,\,$ if 
$i \ne j$ are members of $W$ then $j \notin u_i$.
\medskip
\noindent
[Why?  By de Bruijn and Erd\"os \cite{BrEr}; however we shall
give a proof when we weaken the bound.  First weaken the demand to 
\medskip
\roster
\item "{$(iii)'$}"  $i \in W \and j \in W \and i < j \Rightarrow 
j \notin u_i$. \newline
This we get as follows: choose the $i$-th member by induction.  Next we 
find $W' \subseteq W$ such that $W'$ satisfies (iii); then choose this 
is done similarly
but we choose the members from the top down (inside $W$) so the requirement
on $i$ is $i \in W \and (\forall j)(i < j \in W' \rightarrow i \notin u_j)$
so our situation is similar.  So we have proved the existence, except that 
we get a somewhat weaker bound, which is immaterial here].
\endroster
\medskip

\noindent
Now for some $W' \subseteq W$
\medskip
\roster
\item "{$(*)$}"  $W' \subseteq W,|W'| \ge \frac12 |W|$, and all the
$G \restriction B_i$ for $i \in W'$ are complete graphs or all are
independent sets.
\endroster
\medskip

\noindent
By symmetry we may assume the former. \newline
Let us sum up the relevant points:
\medskip
\roster
\item "{$(A)$}"  $W' \subseteq \{0,\dotsc,i^* - 1\}$, \newline
$|W'| \ge \frac{(m^*_2c_5 - \frac{c_3}{\text{log }n}) \cdot n}
{2(m^*_1)^2(m^*_1 -1)c_4(\text{log }n)^2}$
\smallskip
\item "{$(B)$}"  $G \restriction B_i$ is a complete graph for $i \in W'$
\smallskip
\item "{$(C)$}"  $B_i = \{x_{i,\ell}:\ell < m^*_1\}$ without repetition
and \newline
$i_1,i_2 < i^*,\ell_1,\ell_2 < m^*_1 \Rightarrow x_{i_1,\ell_1} E^G 
x_{i_2,\ell_2} \equiv x_{i_1,0} E^G x_{i_2,0}$.
\endroster
\medskip

\noindent
But by the choice of $m^*_1$ (and as $n$ is large enough hence $|W'|$ 
is large enough) we know 
$|W'| \rightarrow \left( \frac{c_1}{m^*_1} \text{ log }n,\frac{c_1}{1}
\text{ log } n \right)^2$. \newline
We apply it to the graph $\{x_{i,0}:i \in W'\}$. \newline
So one of the following occurs:
\medskip
\roster
\item "{$(\alpha)$}"  there is $W'' \subseteq W'$ such that 
$|W''| \ge \frac{c_1}{m^*_1}
\text{ log }n$ and $\{x_{i,0}:i \in W''\}$ is a complete graph
\endroster
\noindent
or
\roster
\item "{$(\beta)$}"  there is $W'' \subseteq W'$ such that $|W'| \ge c_1
(\text{log }n)$ and $\{x_{i,0}:i \in W''\}$ is a graph with no edges.
\endroster
\medskip

Now if possibility $(\beta)$ holds, then $\{x_{i,0}:i \in W''\}$ is as
required and if possibility $(\alpha)$ holds then $\{x_{i,t}:i \in W'',
t < m^*_1\}$ is as required (see (C) above).
\enddemo
\bigskip
\newpage
\newpage
    
REFERENCES.  
\bibliographystyle{lit-plain}
\bibliography{lista,listb,listx,listf,liste}

\enddocument

\bye

\newpage
    
REFERENCES.  
\bibliographystyle{lit-plain}
\bibliography{lista,listb,listx,listf,liste}

\enddocument

\bye